\documentclass{article}

\usepackage[english]{babel}

\usepackage[letterpaper,top=2cm,bottom=2cm,left=3cm,right=3cm,marginparwidth=1.75cm]{geometry}

\usepackage{amsmath,amsfonts,amssymb, stmaryrd,mathtools,amsthm}
\usepackage{graphicx}
\usepackage[colorlinks=true, allcolors=blue]{hyperref}
\allowdisplaybreaks[4]

\newtheorem{theorem}{Theorem}[section]

\newtheorem{proposition}[theorem]{Proposition}

\newtheorem{conjecture}[theorem]{Conjecture}

\newcommand{\n}{\nabla}

\newcommand{\D}{\overline{\nabla}}

\newcommand{\R}{\mathbb{R}}

\title{Monotonicity of Perelman $\mathcal{W}$-Entropy of Mean Curvature Flow}
\author{ Xiang-Dong Li \footnote{1. State Key Laboratory of Mathematical Sciences, Academy of Mathematics
and Systems Science, Chinese Academy of Sciences, No. 55, Zhongguancun East Road, Beijing,
100190, China; 2. School of Mathematical Sciences, University of Chinese Academy of Sciences,
Beijing, 100049, China. Research supported by National Key R\&D Program of China (No. 2020YF0712702), NSFC No.
 12171458, and Key Laboratory RCSDS, CAS, No. 2008DP173182.  \texttt{xdli@amt.ac.cn}}  \quad and \quad Qi Yan \footnote{Academy of Mathematics
and Systems Science, Chinese Academy of Sciences, No. 55, Zhongguancun East Road, Beijing,
100190.  \texttt{yanqi19@mails.ucas.ac.cn}}     }

\begin{document}
\maketitle

\begin{abstract}
In this paper, we study Perelman' s $ \mathcal{W}$ entropy for mean curvature flow in $\mathbb{R}^{n+1}$. 
Analogously to Perelman's $\mathcal{W}$-entropy  defined for Ricci flow, K. Ecker in \cite{Ecker07} defined a functional $\mathcal{W}$ for the mean curvature flow in $\mathbb{R}^{n+1}$ and the region it encloses, and made the conjecture that this functional is monotonically increasing in time. We modify K. Ecker's definition and, using Hamilton's Harnack inequality for mean curvature flow, prove that our redefined $\mathcal{W}$-entropy is monotonically decreasing in time. Additionally, we provide a rigidity theorem for this $\mathcal{W}$-entropy.
\\\\
{\bf Ketwords:} mean curvature flow, Perelman $\mathcal{W}$-entropy, monotonicity formula, Harnack inequality
\end{abstract}

\section{Introduction}
Let $(M,g)$ be an $n$-dimensional compact Riemannian manifold, $\mathrm{Ric}$ and $R$ be the Ricci curvature and scalar curvature of the Riemannian metric $g$ respectively. Perelman in his 2002 preprint \cite{Perelman02} introduced the following functional called the $\mathcal{W}$ entropy for the Ricci flow
\[
\partial _tg=-2\mathrm{Ric},
\] 
\[
\mathcal{W}(g, f, \tau)=\int_M\left[\tau\left(|\nabla f|^2+R\right)+f-n\right] u d V,
\]
where $\tau>0$ is a time parameter, $dV$ is the volume measure on $(M,g)$, $ f \in C^{\infty}(M)$ satisfies $u=(4\pi\tau)^{-n/2}e^{-f}$ and \[\int_M u dV=1.\]

Perelman proved: if $(g(t),f(t),\tau(t))$ is a solution to the system
\[
\partial_tg = -2\mathrm{Ric},\quad \partial_tf=-\Delta f+|\nabla f|^2-R+\frac{n}{2\tau},\quad \partial_t\tau=-1.
\] 
then the following Perelman $\mathcal{W}$-entropy formula holds
\begin{equation}
    \frac{d}{dt}\mathcal{W}(g(t),f(t),\tau(t))=\int_M2\tau\left|\mathrm{Ric}+\nabla^2f-\frac{g}{2\tau}\right|^2udV.
\end{equation}
This shows that the $\mathcal{W}$-entropy is monotonically increasing along the Ricci flow with respect to $\tau$, and strictly monotonically increasing unless $M$ is a shrinking Ricci soliton:
\[
\mathrm{Ric}+\nabla^2f=\frac{g}{2\tau}.
\]
Perelman used this entropy formula to prove the local non-collapsing theorem for Ricci flow \cite[Section 4]{Perelman02}, which played a crucial role in the final resolution of the Poincar\'{e} conjecture.

Following Perelman, Lei Ni in 2004 studied the $\mathcal{W}$-entropy for the linear heat equation on an $n$-dimensional complete Riemannian manifold $(M,g)$ in \cite{Ni04}. He considered the linear heat equation
\[
\partial_tu=\Delta u
\] 
with positive solution satisfying $\int_M u(x,0)dV=1$ 
\[
u(x,\tau)=\frac{e^{-f(x,\tau)}}{(4\pi\tau)^{n/2}}.
\]
and defined the $\mathcal{W}$ entropy as
\[
\mathcal{W}(f,\tau)=\int_M\left(\tau|\nabla f|^2+f-n\right)udV,
\]
where $\tau>0$ and $d\tau/dt=1$. He proved the following formula
\[
\frac{d}{dt}\mathcal{W}(f,\tau)=-2\int_M\tau\left(\left|\nabla^2f-\frac{g}{2\tau}\right|^2+\mathrm{Ric}(\nabla f,\nabla f)\right)udV.
\]
This implies that on a complete Riemannian manifold with non-negative Ricci curvature, the $\mathcal{W}$-entropy for the linear heat equation is monotonically decreasing with respect to $\tau$.

Inspired by the work of Perelman and Lei Ni, Xiang-Dong Li in \cite{Li12} studied the $\mathcal{W}$-entropy for the weighted heat equation
\[
\partial_tu=Lu
\]
on an $n$-dimensional Riemannian manifold $(M,g)$, where the operator $L=\Delta-\nabla f\cdot\nabla$ is the weighted Laplace operator corresponding to the weighted measure $d\mu=e^{-f(x)}dV$, with $f\in C^{\infty}(M)$. The $m$-dimensional Bakry-\'Emery curvature of $L$ is defined as:
\[
\mathrm{Ric}_{m,n}(L):=\mathrm{Ric}+\nabla^2f-\frac{\nabla f\otimes\nabla f}{m-n}.
\]
He considered the positive solution of the weighted heat equation satisfying $\int_Mu(x,t)d\mu=1$
\[
u(x,t)=\frac{e^{-\phi(x,t)}}{(4\pi t)^{m/2}}
\]
and defined the $\mathcal{W}$ entropy
\[
\mathcal{W}_m(\phi,t)=\int_M\left(t|\nabla \phi|^2+\phi-m\right)\frac{e^{-\phi}}{(4\pi t)^{m/2}}d\mu.
\]
He proved:
\begin{equation}
    \begin{aligned}
        \frac{d}{dt}\mathcal{W}_m(\phi,t)=&-2\int_Mt\left(\left|\nabla^2\phi-\frac{g}{2t}\right|^2+\mathrm{Ric}_{m,n}(L)(\nabla\phi,\nabla\phi)\right)ud\mu\\
        &-\frac{2}{m-n}\int_Mt\left(\nabla f\cdot\nabla\phi+\frac{m-n}{2t}\right)^2ud\mu.
    \end{aligned}
\end{equation}
This implies that on a complete Riemannian manifold with non-negative $m$-dimensional Bakry-\'Emery curvature $\mathrm{Ric}_{m,n}(L)$, the $\mathcal{W}_m$-entropy for the weighted linear heat equation is monotonically decreasing in time. Moreover, a rigidity theorem was obtained: $\frac{d}{dt}\mathcal{W}_m(u,t)=0$ if and only if $M=\R^n$, $m=n$, $f=C$ is constant, and $u$ is the heat kernel.

K. Ecker in \cite{Ecker07} studied the $\mathcal{W}$-entropy for evolving domains in $\mathbb{R}^{n+1}$. Specifically, he considered a family of bounded open subsets $(\Omega_t)_{t\in[0,T]} \subset \mathbb{R}^{n+1}$ whose smooth boundary hypersurfaces $(M_t)_{t\in[0,T]}$ evolve with normal speed $\beta=\beta_{M_t}$:
\[
\beta_{M_t}=-\frac{\partial x}{\partial t}\cdot N,
\]
where $x$ is the embedding of $M_t$, $N$ is the outward unit normal to $\Omega_t$. He defined
\begin{equation}\label{W-entropy}
\mathcal{W}_\beta(\Omega_t, f(t), \tau(t))=\int_{\Omega_t}\left(\tau|\D f|^2+f-(n+1)\right) u d V+2 \tau \int_{M_t} \beta u d S,
\end{equation}
where $\D$ denotes the gradient in $\R^{n+1}$, \[
u=\frac{e^{-f}}{(4\pi\tau)^{(n+1)/2}}.
\]
Assuming $\tau(t)>0$ satisfies $\frac{\partial\tau}{\partial t}=-1$, $f$ satisfies the evolution equation in $\Omega_t$
\[
\frac{\partial f}{\partial t}+\overline{\Delta}f=|\D f|^2+\frac{n+1}{2\tau}.
\] 
Here $\overline{\Delta}$ denotes the Laplace operator in $\R^{n+1}$,
with Neumann boundary condition on $M_t$
\[
\D f\cdot N=\beta,
\] 
and let 
\[
\frac{\partial x}{\partial t}=-\D f(x,t).
\] 

Then K. Ecker obtained
\begin{equation}
\begin{aligned}
\frac{d}{dt}\mathcal{W}_{\beta}(\Omega_t,f(t),\tau(t))
&=2\tau\int_{\Omega_t}\left|\D^2f-\frac{\operatorname{Id}}{2\tau}\right|^2udV\\
&\qquad+2\tau\int_{M_t}\left(\frac{\partial\beta}{\partial t}-2\nabla\beta\cdot\nabla f+h(\nabla f,\nabla f)-\frac{\beta}{2\tau}\right)udS.
\end{aligned}
\end{equation}
where $\nabla$ denotes the covariant derivative on the hypersurface $M_t$.
In the important case of mean curvature flow (i.e., $\beta=H$, the mean curvature of $M_t$), the Harnack quadratic
\[
Z(-\nabla f)=\frac{\partial H}{\partial t}-2\nabla H\cdot\nabla f+h(\nabla f,\nabla f)
\] 
is a key quantity in Hamilton's Harnack inequality for mean curvature flow. 

Hamilton's Harnack inequality for mean curvature flow is a Li-Yau type Harnack inequality, which states that, any weakly convex solution $M_t$ of mean curvature flow 
\[
\frac{\partial x}{\partial t}=-HN,
\] satisfies that for all $t>0$ and all tangent vector $V$ on $M_t$, the following inequality holds:
\begin{equation}\label{Harnack inequality for mean curvature flow}
    \partial_t H+2DH(V)+h(V,V)+\frac{H}{2t}\geq 0.
\end{equation}
where $h$ is the second fundamental form of $M_t$. For more details, we refer to \cite{Hamilton95}

K. Ecker proposed the following conjecture in \cite{Ecker07}:
\begin{conjecture}
In the case of mean curvature flow for compact embedded hypersurfaces $M_t$ in $\R^{n+1}$, for the above $f(t),\tau(t)$ we have
\begin{equation}
\frac{d}{dt}\mathcal{W}(\Omega_t,f(t),\tau(t))
=2\tau\int_{\Omega_t}\left|\D^2f-\frac{\operatorname{Id}}{2\tau}\right|^2udV+2\tau\int_{M_t}\left(Z(-\nabla f)-\frac{H}{2\tau}\right)udS\geq 0.
\end{equation}
\end{conjecture}

\section{Main Results}
We consider a family of bounded open subsets $(\Omega_t)_{t\in [0,T]}\subset\R^{n+1}$ whose smooth boundary hypersurfaces $(M_t)_{t\in [0,T]}$ evolve with normal speed $\beta=\beta_{M_t}$:
\[
\beta_{M_t}=-\frac{\partial x}{\partial t}\cdot N,
\]
where $x$ is the embedding of $M_t$, $N$ is the outward unit normal to $\Omega_t$. 

Different from K. Ecker's definition, we define the $\mathcal{W}$-entropy by:
\begin{equation}
\mathcal{W}_{\beta}(\Omega_t,f(t),t)=\int_{\Omega_t}\left(t|\D f|^2+f-(n+1)\right)udV-2t\int_{M_t}\beta udS.
 \end{equation}
where 
\[
u=\dfrac{e^{-f}}{(4\pi t)^{(n+1)/2}}
\]is a positive solution of the heat equation on $\R^{n+1}$ 
\[
\partial_tu=\overline{\Delta} u.
\] 
Then $f$ satisfies
 \[
\partial_tf-\overline{\Delta}f+|\D f|^2+\frac{n+1}{2t}=0,
\]
We let $\Omega_t$ evolve along the grandient flow of $f$, i.e.
\[
\frac{\partial x}{\partial t}=\D f(x,t),\quad  x\in \Omega_t,
\]
the evolution of the boundary $M_t$ should be compatible with the interior, hence, the normal speed $\beta_{M_t}$ should satisfies the  Neumann boundary condition 
\[
\beta=-\D f\cdot N.
\]

We derive
\begin{equation}
\begin{aligned}
&\frac{d}{dt}\mathcal{W}_{\beta}(\Omega_t,f(t),t)\\
&=-2t\int_{\Omega_t}\left|\D^2f-\frac{\operatorname{Id}}{2t}\right|^2udV\\
&\qquad-2t\int_{M_t}\left(\partial_t\beta+2\nabla\beta\cdot\nabla f+h(\nabla f,\nabla f)+\frac{\beta}{2t}\right)udS.
\end{aligned}
\end{equation}

For mean curvature flow case (i.e. $\beta=H$), we get:
\begin{theorem}[Main Theorem]
    If the compact embedded hypersurface $M_t$ is a weakly convex (i.e., $H\geq 0$) mean curvature flow, with the above setting, we have 
\[
\frac{d}{dt}\mathcal{W}_H(\Omega_t,f(t),t)=-2t\int_{\Omega_t}\left|\D^2f-\frac{\operatorname{Id}}{2t}\right|^2udV-2t\int_{M_t}\left(Z(\nabla f)+\frac{H}{2t}\right)udS\leq 0.
\]
\end{theorem}

\begin{theorem}[Rigidity Theorem]
\[
\frac{d}{dt}\mathcal{W}_H(\Omega_t,f(t),t)=0
\]
if and only if the hypersurface $M_t=\partial\Omega_t$ is the $n$-dimensional sphere $\mathbb{S}(\sqrt{2nt})$ of radius $\sqrt{2nt}$, and
\[
u=\frac{e^{-f}}{(4\pi t)^{\frac{n+1}{2}}}=\frac{e^{-\frac{|x|^2}{4t}}}{(4\pi t)^{\frac{n+1}{2}}}
\]
is the heat kernel for the heat equation on $\R^{n+1}$: \[
\partial_tu=\overline{\Delta}u
\]

\end{theorem}
 We give an intuitive explanation for why we redefine $\mathcal{W}$. Entropy, derived  from thermodynamics and statistical mechanics, is a quantity that describes the degree of disorder in a system. The more chaotic the system , the greater the entropy, the more ordered the system, the smaller the entropy. Huisken's famous work in \cite{Huisken84} states that mean curvature flow with a  initial convex hypersurface will shrink to a circular point in  finite time. This process can be understood as a transition from relative disorder (any initial convex hypersurface) to order (the point at which the mean curvature shrinks). It seems that the entropy of the system should be monotonically decreasing. Since K. Ecker tried to adopt the backward heat equation to prove his $\mathcal{W}$-entropy is monotonically increasing. We just consider the forward heat equation.

\section{Proof of the Main Result}

\begin{proof}
First, we define
\begin{equation}
    W(f):=t(2\overline{\Delta} f-|\D f|^2)+f-(n+1).
\end{equation}
If $\beta$ satisfies
\[
\beta=-\D f\cdot N,
\]
then by the divergence theorem, we have
\begin{equation}
    \mathcal{W}_{\beta}(\Omega_t,f(t),t)=\int_{\Omega_t}WudV.
\end{equation}
We compute
\begin{align*}
    (\partial_t-\overline{\Delta})\overline{\Delta}f&=\partial_t\overline{\Delta}f-\overline{\Delta}^2f =\overline{\Delta}(\partial_t\overline{\Delta})f=-\overline{\Delta}|\D f|^2,\\
    (\partial_t-\overline{\Delta})|\D f|^2&=2\D f\cdot\D\partial_tf-\overline{\Delta}|\D f|^2=2\D f\cdot\D(\overline{\Delta}f-|\D f|^2)-\overline{\Delta}|\D f|^2.
\end{align*}
Thus we have
\begin{align*}
	(\partial_t-\overline{\Delta})(Wu)&=u(\partial_t-\overline{\Delta})W-2\D u\cdot\D W\\
	&=u(\partial_t-\Delta)\left(t(2\overline{\Delta}f-|\D f|^2)+f\right)+2u\D f \cdot\left(2t\D \overline{\Delta}f-t\D|\D f|^2+\D f\right)\\
	&=u\left(2\overline{\Delta}f-|\D f|^2-t\overline{\Delta}|\D f|^2-2t\D f\cdot\D \overline{\Delta}f+2t\D f\cdot\D|\D f|^2-|\D f|^2-\frac{n+1}{2t}\right)\\
	&\qquad+2u\left(2t\D f\cdot\D Lf-t\D f\cdot\D|\D f|^2+|\D f|^2\right)\\
	&=u\left(2\overline{\Delta}f-t\overline{\Delta}|\D f|^2+2t\D f\cdot\D \overline{\Delta}f-\frac{n+1}{2t}\right).
\end{align*}
By the Bochner formula:
\[
\frac{1}{2}\overline{\Delta}|\D f|^2=\langle\D f,\D\overline{\Delta}f\rangle+|\D^2f|^2.
\]
We have:
\begin{equation}
    \begin{aligned}
    (\partial_t-\overline{\Delta})(Wu)&=u\left(2\Delta f-2t|\D^2f|^2-\frac{n+1}{2t}\right)\\
    &=-2tu\left(|\D^2f|^2-\frac{\overline{\Delta}f}{t}+\frac{n+1}{4t^2}\right)\\
    &=-2tu\left|\D^2f-\frac{\mathrm{Id}}{2t}\right|^2.
    \end{aligned}
\end{equation}
Assume that inside $\Omega_t$ the flow follows the gradient flow of $f$, i.e.,
\[
\frac{\partial x}{\partial t}=\D f.
\]
where $x$ is the position coordinate.
Then the evolution of the volume element is given by
\[
\frac{d}{dt}dV=\mathrm{div}(\D f)dV=\overline{\Delta}f dV=\left(\frac{|\D u|^2}{u^2}-\frac{\overline{\Delta} u}{u}\right)dV.
\]
We can proceed to compute the time derivative of $\mathcal{W}$:
\begin{align*}
	\frac{d}{dt}\mathcal{W}_{\beta}&(\Omega_t,f_t,t)=\frac{d}{dt}\int_{\Omega_t}WudV\\
	&=\int_{\Omega_t}\frac{d}{dt}(Wu)d\mu+Wu\frac{d}{dt}dV\\
	&=\int_{\Omega_t}\partial_t(Wu)+\D (Wu)\cdot\D f+W\left(\frac{|\D u|^2}{u}-\overline{\Delta} u\right)dV\\
	&=\int_{\Omega_t}(\partial_t-\overline{\Delta})(Wu)+\overline{\Delta}(Wu)+\D W\cdot u\D f+W\D u\cdot \D f+W\frac{|\D u|^2}{u}-W\overline{\Delta} udV\\
	&=\int_{\Omega_t}(\partial_t-\overline{\Delta})(Wu)dV+\int_{\Omega_t}(\overline{\Delta} Wu+\D W\cdot\D u)dV\\
	&=\int_{\Omega_t}(\partial_t-\overline{\Delta})(Wu)dV+\int_{M_t}\D W\cdot NudS.
\end{align*}
Now we compute $\D W\cdot N.$
Since
\[
	\frac{df}{dt}=\partial_tf+|\D f|^2.
\]
We have 
\[
	W=2t\overline{\Delta}f-t|\D f|^2+f+n+1=2t\frac{df}{dt}-t|\D f|^2+f.
\]
Therefore,
\[
	\D W\cdot N=2t\D\frac{df}{dt}-2t\D^2f(\D f,N)+\D f\cdot N.
\]
Thus,
\begin{align*}
	\frac{d}{dt}\D f&=\partial_t\D f+\D^2f(\D f)=\D\partial_tf+\D^2f(\D f)=\D\left(\frac{df}{dt}-|\D f|^2\right)+\D^2f(\D f)\\
	&=\D\frac{df}{dt}-\D^2f(\D f).
\end{align*}
Also, since the hypersurface $M_t$ evolves according to
\[
	\frac{\partial x}{\partial t}=-\beta N+\nabla f.
\]
where $\nabla$ is the covariant derivative on $M$. Therefore, 
\[
	\frac{dN}{dt}=\nabla\beta+h(\nabla f).
\]
\begin{align*}
	\frac{d \beta}{dt}&=-\frac{d}{dt}(\D f\cdot N)=-\frac{d}{dt}\D f\cdot N-\D f\cdot\frac{dN}{dt}\\
	&=-\D\frac{df}{dt}\cdot N+\D^2f(\D f,N)-\D f\cdot\n\beta-h(\n f,\n f).
\end{align*}
Combining the above computations, we obtain
\begin{align*}
	\D W\cdot N&=2t\D\frac{df}{dt}\cdot N-2t\D^2f(\D f,N)+\D f\cdot N\\
	&=-2t\frac{d\beta}{dt}-2t\n f\cdot\n\beta-2th(\n f,\n f)-\beta\\
	&=-2t\partial_t\beta-4t\n f\cdot\n\beta-2th(\n f,\n f)-\beta\\
	&=-2t\left(\partial_t\beta+2\n f\cdot\n\beta+h(\n f,\n f)+\frac{\beta}{2t}\right).
\end{align*}

Finally, we obtain the following $W$-entropy formula:
\begin{equation}
	\begin{aligned}
		\frac{d}{dt}\mathcal{W}_{\beta}(\Omega_t,f(t),t)=&-2t\int_{\Omega_t}\left|\D^2f-\frac{\operatorname{Id}}{2t}\right|^2udV\\
		&\quad-2t\int_{M_t}\left(\partial_t\beta+2\n\beta\cdot\n f+h(\n f,\n f)+\frac{\beta}{2t}\right)udS.
	\end{aligned}
\end{equation}
For weakly convex mean curvature flow, by Hamilton's Harnack inequality \eqref{Harnack inequality for mean curvature flow}, we have
\begin{equation}
	\begin{aligned}
		\frac{d}{dt}\mathcal{W}_{H}(\Omega_t,f(t),t)=&-2t\int_{\Omega_t}\left|\D^2f-\frac{\operatorname{Id}}{2t}\right|^2udV\\
		&\quad-2t\int_{M_t}\left(\partial_tH+2\n H\cdot\n f+h(\n f,\n f)+\frac{H}{2t}\right)udS\\
		&\leq 0.
	\end{aligned}
\end{equation}

\end{proof}

We define the entropy:
\begin{equation}
    \mu_{\beta}(\Omega_t,t):=\inf_{\int_{\Omega_t}udV=1}\mathcal{W}_{\beta}(\Omega_t,f(t),t).
\end{equation}

Using the method of Lagrange multipliers, we can prove the following:
\begin{proposition}
A function $f:\Omega\to\R$ is a minimizer of $\mu_{\beta}(\Omega_t,t)$ if and only if it satisfies the following conditions: in $\Omega_t$
\begin{equation}
    W(f):=t(2\overline{\Delta} f-|\nabla f|^2)+f-(n+1)=\mu_{\beta}(\Omega_t,t).
\end{equation}
On $M_t=\partial\Omega_t$
\[
\beta=-\nabla f\cdot N.
\]
and
\[
\int_{\Omega_t}\frac{e^{-f}}{(4\pi t)^{(n+1)/2}}=1.
\]
\end{proposition}
For details, see \cite[Proposition 2.1]{Ecker14}.

The Harnack inequality for mean curvature flow
\begin{equation}
    \partial_tH+2\nabla H\cdot V+h(V,V)+\frac{H}{2t}\geq 0.
\end{equation}
holds with equality if and only if $M_t$ is a convex expanding self-similar solution and
\[
V=\frac{x^{\top}}{2t},
\]
where $x$ denotes the position vector in $\R^{n+1}$, $x^{\top}$ denotes its tangential component on the surface $M$.

For the Harnack quantity appearing in our $\frac{d}{dt}\mathcal{W}$
\[
\partial_tH+2\nabla H\cdot\nabla f+h(\nabla f,\nabla f)+\frac{H}{2t},
\]
if it is identically zero, then the hypersurface $M_t$ must be a convex expanding self-similar solution of mean curvature flow, and $\nabla f=\dfrac{x^{\top}}{2t}$. This implies
\[
f =\frac{|x|^2}{4t},
\]
at which point,
\[
u=\frac{e^{-f}}{(4\pi t)^{\frac{n+1}{2}}}=\frac{e^{-\frac{|x|^2}{4t}}}{(4\pi t)^{\frac{n+1}{2}}}
\]
is exactly the heat kernel for the heat equation on $\R^{n+1}$! Moreover, the compact convex expanding self-similar solution of mean curvature flow can only be the $n$-dimensional sphere $\mathbb{S}(\sqrt{2nt})$ of radius $\sqrt{2nt}$, see \cite[Theorem 4.1]{Huisken90}.

Therefore, we have the following rigidity theorem:
\begin{theorem}[Rigidity Theorem]
\[
\frac{d}{dt}\mathcal{W}_H(\Omega_t,f(t),t)=0
\]
if and only if the hypersurface $M_t=\partial\Omega_t$ is the $n$-dimensional sphere $\mathbb{S}(\sqrt{2nt})$ of radius $\sqrt{2nt}$, and
\[
u=\frac{e^{-f}}{(4\pi t)^{\frac{n+1}{2}}}=\frac{e^{-\frac{|x|^2}{4t}}}{(4\pi t)^{\frac{n+1}{2}}}
\]
is the heat kernel for the heat equation on $\R^{n+1}$: \[
\partial_tu=\overline{\Delta}u
\]

\end{theorem}

\bibliographystyle{alpha}
\bibliography{sample}

\end{document}